\documentclass{amsart}

\usepackage[dvips]{color}
\usepackage{amsfonts}
\usepackage{graphicx}
\usepackage{amsmath}
\usepackage{amssymb}
\usepackage{amsthm}
\usepackage{mathrsfs}

\newtheorem{theorem}{Theorem}[section]
\newtheorem{lemma}[theorem]{Lemma}

\theoremstyle{definition}
\newtheorem{definition}[theorem]{Definition}

 \def\HollowBox #1#2{{\dimen0=#1 \advance\dimen0 by -#2       
       \dimen1=#1 \advance\dimen1 by #2                       
        \vrule height #1 depth #2 width #2                    
        \vrule height 0pt depth #2 width #1                   
        \llap{\vrule height #1 depth -\dimen0 width \dimen1}%
       \hskip -#2                                             
       \vrule height #1 depth #2 width #2}}                   
              
\theoremstyle{remark}

\numberwithin{equation}{section}

\newfam\msbfam
\font\tenmsb=msbm10  \textfont\msbfam=\tenmsb
\font\sevenmsb=msbm7  \scriptfont\msbfam=\sevenmsb
\font\fivemsb=msbm5    \scriptscriptfont\msbfam=\fivemsb
\def\Bbb{\fam\msbfam \tenmsb}

\def\CC{{\Bbb C}}

\newfam\msbbfam
\font\tenmsbb=msbm10  scaled \magstep1 \textfont\msbbfam=\tenmsbb
\font\sevenmsbb=msbm7  scaled \magstep1 \scriptfont\msbbfam=\sevenmsbb
\font\fivemsbb=msbm5    scaled \magstep1 \scriptscriptfont\msbbfam=\fivemsbb

\begin{document}

\begin{center}
\Large \bf Restricted Admissible Limit for Domains of Finite Type
\bigskip \\
\rm Steven G. Krantz and Baili Min
\end{center}
\vspace*{.25in}

\begin{quote}
{\bf Abstract:}  We investigate the boundary behavior of holomorphic functions with respect to a family of curves in 
a domain of finite type. This work is a generalization of \u{C}irka's classical result on the unit 
ball and it supplements the result by Cima and Krantz on the Lindel\"{o}f principle for general domains.
See [KRA2] for some recent developments in this subject.

Our discussion is carried out in $\CC^2$.  
\end{quote}

\parskip=\baselineskip
\section{Background}
The classical Lindel\"{o}f principle says that, in the unit disk, if a bounded holomorphic function has a boundary limit $L$ at $\hbox{e}^{i\theta}$ along the radial approach $r\hbox{e}^{i\theta}, r \to 1-$, then it has non-tangential limit $L$ at $\hbox{e}^{i\theta}$. When it comes to several variables, the story is much more subtle. Much has been investigated for the unit ball $B\subset \CC^n$ and \u{C}irka's generalization in [CIR] is of outstanding significance. His version is to consider the function's boundary behavior along  two types of curves: special and restricted ones. Basically speaking, a special curve in $B$ 
is one that is contained in a paraboloid, and a restricted curve is a special one that also has a non-tangential projection in the complex normal direction.  These definitions can also be found in Rudin's book [RUD], as well as the result given by \u{C}irka which supplemented the classical  Lindel\"{o}f principle as follows.
\begin{theorem}[\u{C}irka]  \sl
Let $S$, the unit sphere, be the boundary of the ball $B$.  Suppose that $f \in H^\infty(B)$, $\zeta \in S$, and $\Gamma_0(t)$ is a special curve with $\lim_{t \to 1}\Gamma_0(t)=\zeta$, and
\begin{equation}
\lim_{t \to 1}f(\Gamma_0(t))=L.
\end{equation}
Then, for any restricted curve $\Gamma(t)$ with $\lim_{t \to 1}\Gamma(t)=\zeta$, we have
\begin{equation}
\lim_{t \to 1}f(\Gamma(t))=L.
\end{equation}
\end{theorem} 

One of the important parts of this theorem is that it is
linked to the admissible convergence of bounded holomorphic
functions on the unit ball because of the fact that any
restricted curve will eventually be contained in an admissible
approach region based at $\zeta$, as discussed in [RUD].
Readers are referred to [KRA1] and [RUD] for more about
admissible approach regions.

There are still many unknowns and developments about this principle for more general domains in $\CC^n$, among which Cima and Krantz formulated another version of the Lindel\"{o}f principle from the point
of view of normal functions in their work [CIK].

In this paper, we discuss the generalization to a domain $\Omega \subset \CC^2$ that is of finite type. We first define the generalized special and restricted curves.

Suppose that $\Gamma=\Gamma(t)$ is a curve in $\Omega$ approaching to $\zeta \in \partial \Omega$ and let $\gamma=\gamma(t)$ be its orthogonal projection to the complex line through a transversal vector to $\partial \Omega$ at $\zeta$.
\begin{definition}  \rm
The curve $\Gamma$ is special if $|\Gamma(t) - \gamma(t)|^m = o(|\zeta-\gamma(t)|)$,
where the number $m$ is the type of the point $\zeta$.
\end{definition}
\begin{definition}   \rm
The curve $\Gamma$ is {\it restricted} if it is special, and $\gamma$ is non-tangential.
\end{definition}
Careful discussion of these matters will be carried out in the next section.

The Lindel\"{o}f principle in terms of these curves is then as follows.  This is the main result of the present paper.
We shall explain the concept of finite type in the next section. 

\begin{theorem}	\label{mainthm}  \sl
Suppose that $\Omega \subset \mathbb{C}^2$ is a domain of finite type and $\zeta \in \partial \Omega$. If $f\in H^\infty(\Omega), \zeta \in \partial \Omega$, and $\Gamma_0$ is a restricted curve that approaches to $\zeta$ along which $f$ converges:
\begin{displaymath}
\lim_{t \to 1}f(\Gamma_0(t))=L \, ,
\end{displaymath}
then $f$ has restricted admissible limit $L$ at $\zeta$, which means, for any restricted curve $\Gamma(t)$ in $\Omega$ with $\lim_{t \to 1-} \Gamma(t)=\zeta$, we have the limit
\begin{displaymath}
\lim_{t \to 1}f(\Gamma(t))=L \, .
\end{displaymath}
\end{theorem}

\section{Special and Restricted Curves}

Let $\Omega$ be a bounded domain in $\CC^2$, and suppose that
$\zeta$ is a boundary point in $\partial \Omega$ with
$\nu$ the unit outward normal vector.

The complex tangent space  to $\partial \Omega$ at $\zeta$ is the largest complex subspace of the real
tangent space $T_\zeta$. The complex normal space is just $\CC \nu$.

Suppose that $\Gamma$ is a curve in $\Omega$ approaching $\zeta$,  $\Gamma: [0,1) \to \Omega$, and 
\begin{equation}
\lim_{t \to 1-}\Gamma(t)=\zeta. 
\end{equation}
Then let $\gamma$ be its orthogonal projection to $\CC \nu$
\begin{equation}
\gamma= \gamma(t)=\langle \Gamma(t), \nu \rangle\nu.
\end{equation}

The definition of the type of $\zeta$ is subtle. Geometrically speaking, the type is the order of contact: if $l$ is a nonsingular complex curve tangent to $\partial \Omega$ at $\zeta$, and for $z \in l$, 
\begin{equation}
\text{dist}(z, \partial \Omega)=\mathcal{O}(|z-\zeta|^\tau),
\end{equation}
and if the number $\tau$ cannot be improved, then we say $\zeta$ has the type $\tau$.  See [KRA1] for more on finite type.
Readers are also referred to [JJK], [BLG], [KRA1], and [NSW] for the discussion of the concept of type. 

In his book [RUD], Rudin discussed the relationship between
restricted curve and the admissible approach regions in the
unit ball. In this section, we also hope to set up a similar
context. Since we are dealing with a weakly pseudoconvex
domain while the unit ball is strongly pseudoconvex, the
definition of the admissible approach region is very
different. In our case, the admissible region is broader in
the complex tangential direction and relies on the type of the
boundary point.  See [EMS] and [NSW] for more about this
definition.

One of the results in the paper [DBF] implies that the shape of the admissible approach region $\mathscr{A}_{\alpha}(1,0)$ is reflected by the order of contact such that
if the dimension in the complex tangential direction is $r$ then the dimension in the complex normal direction is $r^m$.

Let us review the definitions of special and restricted curves:
\begin{definition}  \rm
A $\zeta$-curve $\Gamma=\Gamma(t)$ is special if
\begin{equation}
\label{specialcurve}
|\Gamma(t) - \gamma(t)|^m = o(|\zeta-\gamma(t)|).
\end{equation}
\end{definition}

Since $\gamma$ is the projection of $\Gamma$ onto the complex normal direction $\nu$, $\zeta-\gamma(t)$ is exactly the complex normal component at $\zeta$ of $\zeta-\Gamma(t)$, while $\Gamma(t) - \gamma(t)$ is the complex tangential component. The equation ~\ref{specialcurve} relates the complex tangential and normal components. Recall that the type of $\zeta$ also reflects the order of contact in the sense of complex tangential and normal components in a similar way.

However, a special curve may not be in an admissible approach region with vertex $\zeta$, which is tangent to $\partial \Omega$ only in the complex tangential directions.

\begin{definition}   \rm
The curve $\Gamma$ is {\it restricted} if it is special, and $\gamma$ is non-tangential to $\partial \Omega$ at $\zeta$
\end{definition}

With this extra requirement, we have a quick result:
\begin{lemma}	  \sl
\label{curveregion}
If $\Gamma(t)$ is a restricted curve then, for $t$ close enough to 1, $\Gamma(t)$ lies in an admissible approach region $\mathscr{A}_{\alpha}(1,0)$.
\end{lemma} 

Consequently we define
\begin{definition}
A function $f: \Omega \to \mathbb{C}$ has a restricted admissible limit $L$ at $\zeta \in \partial \Omega$ if 
\begin{displaymath}
\lim_{t \to 1}f(\Gamma(t))=L
\end{displaymath}
for every restricted $\zeta$-curve $\Gamma(t)$.
\end{definition}
We observe that admissible convergence implies the restricted admissible convergence, but not vice versa. 

\section{proof of the main theorem}

We now prove our main result, Theorem \ref{mainthm}, with two major steps. 

First we want to see that
\begin{displaymath}
\lim_{t \to 1} \biggl  (f\big(\Gamma(t)\big)-f\big(\gamma(t)\big) \biggr ) = 0 \, .
\end{displaymath}
for \textit{any} special $\zeta$-curve $\Gamma(t)$.

Suppose that $\Gamma(t)$ is a curve in $\Omega$ that approaches $\zeta \in \partial \Omega$ with $\gamma(t)$ as its projection onto $\CC \nu$, the complex normal direction to $\partial \Omega$ at $\zeta$. Fix at $t$  the point $(1-\lambda)\gamma(t)+\lambda\Gamma(t)=\gamma(t)+\lambda(\Gamma(t)-\gamma(t))$, which has the component $\lambda(\Gamma(t)-\gamma(t))$ in the complex tangential direction from $\zeta$ and complex normal component $\zeta-\gamma(t)$.

The type of $\zeta$ is $m$, which is the greatest order of contact to the bondary along a nonsingular complex curve that is tangential to the boundary. Notice that $\Omega \subset \CC^2$ and there is only one complex tangential direction to $\partial \Omega$ at $\zeta$, so we can always do a holomorphic change of coordinates to change the
curve of best contact into the complex line passing through $\zeta$, along the complex tangential direction. Details can be found in [GUN]. Due to this geometry, there exists a constant $k$ such that, if
\begin{equation}
|\lambda|^m|\Gamma(t)-\gamma(t)|^m <k|\zeta-\gamma(t)|,
\end{equation}
then the point $(1-\lambda)\gamma(t)+\lambda\Gamma(t)=\gamma(t)+\lambda(\Gamma(t)-\gamma(t))$ is in $\Omega$. 

Set 
\begin{equation}
R(t)=\frac{(k|\zeta-\gamma(t)|)^{\frac{1}{m}}}{|\Gamma(t)-\gamma(t)|}.
\end{equation}
The argument above indicates that, if $|\lambda|<R$, then $(1-\lambda)\gamma+\lambda\Gamma \in \Omega$.

Since $\Gamma$ is special, we know from equation (\ref{specialcurve}) that 
\begin{equation}
\lim_{t \to 1}\frac{1}{R(t)}=0.
\end{equation}

Now, for each $t$, we can define a holomorphic function
\begin{displaymath}
g(\lambda)=f\big((1-\lambda)\gamma(t)+\lambda\Gamma(t)\big)
\end{displaymath}
and consider the disc $\mathbb{D}_{R(t)}=\{\lambda \in \mathbb{C}: |\lambda|<R(t)\}$. Applying the Schwarz lemma to $g(\lambda)-g(0)$ in $\mathbb{D}_R$, we are able to get
\begin{equation}
|g(1)-g(0)|\leqslant\frac{2\left\|f\right\|_\infty}{R(t)}.
\end{equation}
Then it is obvious that
\begin{equation}
\lim_{t \to 1} (f(\Gamma(t))-f(\gamma(t)))=0.
\end{equation}

Secondly, recall that a restricted curve is special, therefore we can apply the previous result to the restricted curve $\Gamma_0(t)$ given in the hypothesis:
\begin{equation}
\lim_{t \to 1} (f(\Gamma_0(t))-f(\gamma_0(t)))=0.
\end{equation}
By the hypothesis of the theorem,
\begin{equation}
\lim_{t \to 1}f(\Gamma_0(t))=L,
\end{equation}
we immediately know that
\begin{equation}
\lim_{t \to 1}f(\gamma_0(t))=L.
\end{equation}

Now suppose that $\Gamma$ is an arbitrary restricted curve approaching to $\zeta$. 
Notice that both $\gamma$ and $\gamma_0$ are non-tangential at $\zeta$ according to the definition of being restricted, so by the Lindel\"{o}f principle (see [HIL]) we know that
\begin{equation}
\lim_{t \to 1}f(\gamma(t))=\lim_{t \to 1}f(\gamma_0(t))=L,
\end{equation}
and it then follows from the equation that
\begin{equation}
\lim_{t \to 1} f(\Gamma(t))=L,
\end{equation}
That concludes the proof.

\section{Concluding Remarks}

We have studied in this paper the detailed behavior of a holomorphic function
on a special curve and a restricted curve.  We have restricted attention
to finite type domains in $\CC^2$.  There is definite interest in developing
these ideas in higher dimensions, and on more general types of domains.
\bigskip \bigskip \bigskip \\

\noindent {\large \sc References}
\vspace*{.2in}

\begin{enumerate}

\item[{\bf [BAS]}]  F. Bagemihl and W. Seidel, Some boundary
properties of analytic functions, {\it Math.\ Z.} 61(1954), 186--199.

\item[{\bf [BLG]}] T. Bloom and I. Graham, A geometric characterization of points of type  $m$ on real submanifolds of $\CC^n$,  {\it Jour. Diff. Geom.} 45(1978), 133--147

\item[{\bf [CIR]}] E.-M. \u{C}irka, The Lindel\"{o}f and Fatou theorems in $\CC^n$, {\it
Mat. U.S.S.R. Sb} 92(1973), 622--644. {\it
Math. U.S.S.R. Sb} 21(1973), 619--641.

\item[{\bf [CIK]}] J. C. Cima and S. Krantz, The Lindel\"{o}f principle and normal functions of several complex variables, {\it Duke Math. \ J.} 50(1983), 303--328

\item[{\bf [DBF]}] F. Di Biase and B. Fischer, Boundary behaviour of $H^p$ functions on convex domains of finite type in $\CC^n$, {\it Pacific Jour. Math.} 183(1998), 25--38

\item[{\bf [DOG]}] F. Docquier and H. Grauert, Leisches Problem
und Rungescher Satz f\"{u}r Teilgebiete Steinscher
Mannigfaltigkeiten (German), {\it Math.\ Ann.} 140(1960),
94--123.

\item[{\bf [GUN]}] R. C. Gunning, {\it Lectures on Complex Analytic Varieties: The Local Parame-\\trization Theorem},  Princeton University Press, Princeton, NJ, 1970

\item[{\bf [HIL]}]  E. Hille, {\it Analytic Function Theory}, Ginn, Boston, 1959.

\item[{\bf [JJK]}] J. J. Kohn, Boundary behavior of $\overline{\partial}$ on weakly pseudo-convex manifolds of dimension two,  {\it Jour. Diff. Geom.} 6(1972), 523--542

\item[{\bf [KRA1]}] S. G. Krantz, {\it Function Theory of
Several Complex Variables}, 2nd ed., American Mathematical
Society, Providence, RI, 2001.

\item[{\bf [KRA2]}] S. G. Krantz, The Lindel\"{o}f principle in
several complex variables {\it J. Math. Anal. Appl.}
326(2007), 1190--1198

\item[{\bf [NAG]}] A. Nagel, Smooth zero sets and interpolation
sets for some algebras of holomorphic functions on strictly
pseudoconvex domains, {\it Duke Math.\ J.} 43(1976), 323--348.

\item[{\bf [NSW]}] A. Nagel, E. M. Stein, and S. Wainger, Boundary behavior of functions holomorphic on domains of finite type, {\it Proc. Nat. Acad. Sci.} 78(1981), no. 11, part 1, 6596--6599.

\item[{\bf [RUD]}] W. Rudin, {\it Function Theory in the Unit
Ball of $\CC^n$}, reprint of the 1980 edition, Classics in Mathematics,
Springer-Verlag, Berlin, 2008.

\item[{\bf [EMS]}] E. M. Stein, {\it Boundary Behavior of Holomorphic Functions of Several Complex Variables}, Princeton University Press, Princeton, 1970.

\end{enumerate}
\vspace*{.25in}

\begin{quote}
S. G. Krantz  \\
Department of Mathematics \\
Washington University in St. Louis \\
St.\ Louis, Missouri 63130  \\
{\tt sk@math.wustl.edu}
\end{quote}

\begin{quote}
Baili Min \\
Department of Mathematics \\
Lafayette College \\
Easton, PA 18042 \\
{\tt minbaili@gmail.com}
\end{quote}

\end{document}